\documentclass{amsart}
\usepackage{amsthm,amsmath,amsfonts,amssymb,amscd,mathrsfs,graphics}

\usepackage{latexsym}
\usepackage{graphicx}

\usepackage{caption}

\usepackage{epic}

\usepackage{enumerate}
\usepackage[...]{youngtab}

\usepackage{placeins}
\usepackage{tikz}
\usepackage{listings}

\usepackage{todonotes}

\newtheorem{thm}{Theorem}[section]

\newtheorem{lemma}[thm]{Lemma}
\newtheorem{prop}[thm]{Proposition}
\newtheorem{cor}[thm]{Corollary}

\makeatletter
\renewcommand{\@seccntformat}[1]{\S{\csname
the#1\endcsname}\hspace{0.5em}}
\makeatother

\begin{document}

\title{Strong Gelfand Pairs of the Dihedral and Dicyclic Groups}

\author{Joseph E. Marrow}
  \address{Department of Mathematics,  Brigham Young University, Provo, 
UT 84602, U.S.A.
E-mail:jemarrow@math.byu.edu}
\date{}
\maketitle

\begin{abstract}  
A strong Gelfand pair $(G, H)$ is a finite group $G$ and a subgroup $H$ where every irreducible character of $H$ induces to a multiplicity-free character of $G$. We determine the strong Gelfand pairs of the dihedral groups $\mathrm{D}_{2n}$ and the dicyclic groups $\mathrm{Dic}_{4n}$ for all $n$.

\medskip

\medskip
\end{abstract}

\theoremstyle{plain}

\theoremstyle{definition}
\newtheorem*{dfn}{Definition}
\newtheorem{exa}[thm]{Example}
\newtheorem{rem}[thm]{Remark}

\newcommand{\ds}{\displaystyle}
\newcommand{\bs}{\boldsymbol}
\newcommand{\mb}{\mathbb}
\newcommand{\mc}{\mathcal}
\newcommand{\mf}{\mathfrak}
\renewcommand{\mod}{\operatorname{mod}}
\newcommand{\mult}{\operatorname{Mult}}

\def \a{\alpha} \def \b{\beta} \def \d{\delta} \def \e{\varepsilon} \def \g{\gamma} \def \k{\kappa} \def \l{\lambda} \def \s{\sigma} \def \t{\theta} \def \z{\zeta}

\numberwithin{equation}{section}

\setlength{\leftmargini}{1.em} \setlength{\leftmarginii}{1.em}
\renewcommand{\labelenumi}{\setlength{\labelwidth}{\leftmargin}
   \addtolength{\labelwidth}{-\labelsep}
   \hbox to \labelwidth{\theenumi.\hfill}}

\maketitle

\section{Introduction}
A \textit{Gelfand pair} $(G, H)$ is a finite group $G$ and a subgroup $H \leq G$ where the trivial character of $H$ induces to a multiplicity-free character of $G$. There are many conditions equivalent to being a Gelfand pair; see \cite{grady, hu, harm}. This paper is adapted from the author's Master's thesis \cite{thesis}.

In \cite{hu, GH} every strong Gelfand pair of $\mathrm{SL}_2(p^n)$ is determined, for $p$ a prime. For finite groups $G, H$ with $H\leq G$, the pair $(G, H)$ is called a \textit{strong Gelfand pair} if every irreducible character of $H$ induces a multiplicity-free character of $G$, i.e. for all $\psi \in \hat{H}$ and $\chi\in\hat{G}$, we have $\langle \psi \uparrow G, \chi\rangle \leq 1$. 

Generators $a, b$ for $\mathrm{D}_{2n}$ and $\mathrm{Dic}_{4n}$ are defined in \S $2$. For the dihedral groups $\mathrm{D}_{2n}$ we will prove Theorem \ref{dsgp}.

\begin{thm}\label{dsgp}
For $n\geq 1$, and $H \leq \mathrm{D}_{2n}$, if $(\mathrm{D}_{2n}, H)$ is a strong Gelfand pair, then $H$ is one of \\(i) a reflection subgroup; \\(ii) a dihedral subgroup; \\(iii) the maximal cyclic subgroup; or \\(iv) if $n$ is even the index four cyclic subgroup.
\end{thm}

For the dicyclic groups $\mathrm{Dic}_{4n}$ we show Theorem \ref{dicsgp}.
\begin{thm}\label{dicsgp}
For $n\geq 2$, and $H\leq \mathrm{Dic}_{4n}$ if $(\mathrm{Dic}_{4n}, H)$ is a strong Gelfand pair, then $H$ is one of \\(i) a subgroup of the form $\langle ba^i\rangle$; \\(ii) a dicyclic subgroup; \\(iii) cyclic of order $n$ or $2n$. 
\end{thm}
If $n=1$, then $\mathrm{Dic}_4\cong\mathcal{C}_4$ and all subgroups are strong Gelfand.

\section{Definitions and Preliminary Results}
For the dihedral group $\mathrm{D}_{2n}$we use the presentation
$$
\mathrm{D}_{2n} = \langle a, b \,\vert\, a^n = b^2 = 1, bab=a^{-1}\rangle.
$$
For $\mathrm{Dic}_{4n}$ we use 
$$
\mathrm{Dic}_{4n} = \langle a, b \,\vert \, a^n = b^2, b^4=1, bab^{-1}=a^{-1}\rangle.
$$


We will use $a, b$ as the generators for these groups throughout the paper. 
\begin{lemma}\cite{hu}\label{stack}
Suppose $K \leq H \leq G$ are groups. If $(G, H)$ is not a strong Gelfand pair, then neither is $(G, K)$.
\end{lemma}

We define $H\leq G$ to be a \textit{strong Gelfand subgroup} if and only if $(G, H)$ is a strong Gelfand pair. We need the following facts about the dihedral and dicyclic groups:
\begin{lemma}\label{dsub}
Every subgroup of a dihedral group is either cyclic or dihedral.
\end{lemma}
\begin{lemma}\label{dicsub}
Every subgroup of a dicyclic group is either cyclic or dicyclic.
\end{lemma}

The character table for a cyclic group $\mathcal{C}_n$ generated by $g$ is given in Table \ref{C}, where $\zeta=e^{2\pi i/n}$. Recall that the character table for a direct product of groups is the tensor product of their character tables. This allows us to use Table \ref{C} to obtain the character table for $V_4=\mathcal{C}_2\times\mathcal{C}_2$.

\begin{table}[h!]
\centering
\caption{Character Table for $\mathcal{C}_n$}
\label{C}
\vspace{5pt}
\begin{tabular}{c|cccc}
 & $1$ & $g$ & $\ldots$ & $g^{n-1}$\\
\hline
$\mu_0$ & $1$ & $1$ & $\ldots$ & $1$\\
$\mu_1$ & $1$ & $\zeta$ & $\ldots$ & $\zeta^{(n-1)}$\\
$\mu_2$ & $1$ & $\zeta^2$ & $\ldots$ & $\zeta^{2(n-1)}$\\
$\vdots$ & $\vdots$ & $\vdots$ & $\ddots$ & $\vdots$\\
$\mu_{n-1}$ & $1$ & $\zeta^{n-1}$ & $\ldots$ & $\zeta^{(n-1)^2}$
\end{tabular}
\end{table}

\section{Dihedral Groups}
The character table for the dihedral group $\mathrm{D}_{2n}, n\geq 3,$ is given in \cite{jl} and depends on the parity of $n$. If $n$ is odd the character table is given in Table \ref{dodd}, where $1 \leq r, j \leq (n-1)/2$.

\begin{table}[h!]
\centering
\caption{Character Table for $\mathrm{D}_{2n}$ with $n$ odd}
\label{dodd}
\vspace{5pt}
\begin{tabular}{c|ccc}
 & $1$ & $a^r$ & $b$\\
\hline
$\chi_1$ & $1$ & $1$ & $1$\\
$\chi_2$ & $1$ & $1$ & $-1$\\
$\psi_j$ & $2$ & $e^{\frac{2\pi i r j}{n}}+e^{-\frac{2\pi i r j}{n}}$ & $0$
\end{tabular}
\end{table}

If $n=2m$ the character table is given in Table \ref{deven} where $1 \leq r, j \leq m-1$.

\begin{table}[h!]
\centering
\caption{Character Table for $\mathrm{D}_{2n}$ with $n=2m$}
\label{deven}
\vspace{5pt}
\begin{tabular}{c|ccccc}
 & $1$ & $a^m$ & $a^r$ & $b$ & $ab$\\
\hline
$\chi_1$ & $1$ & $1$ & $1$ & $1$ & $1$\\
$\chi_2$ & $1$ & $1$ & $1$ & $-1$ & $-1$\\
$\chi_3$ & $1$ & $(-1)^m$ & $(-1)^r$ & $1$ & $-1$\\
$\chi_4$ & $1$ & $(-1)^m$ & $(-1)^r$ & $-1$ & $1$\\
$\psi_j$ & $2$ & $2(-1)^j$ & $e^{\frac{2 \pi i rj}{n}}+e^{-\frac{2\pi i rj}{n}}$ & $0$ & $0$
\end{tabular}
\end{table}

When $n=2$, $D_{4}\cong V_4=\{1, a, b, ab\}$ and when $n=1$ we have $D_2=\mathcal{C}_2$. Since these groups are abelian, it is immediate that all subgroups are strong Gelfand pairs. As all proper subgroups of $D_4$ are reflection subgroups, we can limit our considerations to the case when $n\geq 3$. 

\begin{prop}\label{Dsmall}
For $n \geq 3$, the pair $(\mathrm{D}_{2n}, \langle ba^i\rangle)$ is a strong Gelfand pair.
\end{prop}
\noindent\textit{Proof.} 
Suppose first that $n$ is odd; then the character table of $\mathrm{D}_{2n}$ is given in Table \ref{dodd} and all elements $ba^i$ are conjugate to $b$, so without loss of generality we suppose $i=0$. Since $\langle b\rangle\cong \mathcal{C}_2$ we can induce the characters from Table \ref{C} to get Table \ref{cycupdodd}. 

\begin{table}[h!]
\centering
\caption{Characters of $\mathcal{C}_2\cong \langle b\rangle$ induced to $\mathrm{D}_{2n}$ for $n$ odd}
\label{cycupdodd}
\vspace{5pt}
\begin{tabular}{c|ccc}
 & $1$ & $a^r$ & $b$\\
\hline
$\mu_0\uparrow \mathrm{D}_{2n}$ & $n$ & $0$ & $1$\\
$\mu_1\uparrow \mathrm{D}_{2n}$ & $n$ & $0$ & $-1$
\end{tabular}
\end{table}

Then 
$$
\langle \mu_0\uparrow \mathrm{D}_{2n}, \chi_1\rangle=1; \quad
\langle \mu_0\uparrow \mathrm{D}_{2n}, \chi_2\rangle=0; \quad
\langle \mu_0\uparrow \mathrm{D}_{2n}, \psi_j\rangle=1;
$$
\vspace{-12pt}
$$
\langle \mu_1\uparrow \mathrm{D}_{2n}, \chi_1\rangle=0; \quad
\langle \mu_1\uparrow \mathrm{D}_{2n}, \chi_2\rangle=1; \quad
\langle \mu_1\uparrow \mathrm{D}_{2n}, \psi_j\rangle=1;
$$
showing that $(\mathrm{D}_{2n}, \langle ba^i\rangle)$ is a strong Gelfand pair when $n$ is odd. 

Suppose $n=2m$. We wish to show that $(\mathrm{D}_{2n}, \langle ba^i\rangle)$ is a strong Gelfand pair for all $i$. It is no longer the case that all reflections are conjugate, since there are two conjugacy classes of reflections, represented by $b$ and $ab$ in Table \ref{deven}. Without loss of generality we suppose that $i=0, 1$. However, the only characters which disagree on the values of $b$, $ab$ are $\chi_3, \chi_4$. We have $\chi_3(b)=-\chi_4(ab)$ and $\chi_3(ab) = -\chi_4(b)$. Thus it is sufficient to observe merely that $\langle b \rangle \cong \mathcal{C}_2$ and that determining if $(\mathrm{D}_{2n}, \langle b\rangle)$ is a strong Gelfand pair will also determine whether $(\mathrm{D}_{2n}, \langle ab\rangle)$ is a strong Gelfand pair. Inducing the characters of $\langle b \rangle$ to $\mathrm{D}_{2n}$ gives Table \ref{cycupdeven}.

\begin{table}[h!]
\centering
\caption{Characters of $\mathcal{C}_2\cong \langle b\rangle$ induced to $\mathrm{D}_{2n}$ for $n=2m$}
\label{cycupdeven}
\vspace{5pt}
\begin{tabular}{c|ccccc}
 & $1$ & $a^m$ & $a^r$ & $b$ & $ab$\\
\hline
$\mu_0 \uparrow \mathrm{D}_{2n}$ & $n$ & $0$ & $0$ & $2$ & $0$\\
$\mu_1 \uparrow \mathrm{D}_{2n}$ & $n$ & $0$ & $0$ & $-2$ & $0$
\end{tabular}
\end{table}

The inner products of $\mu_0\uparrow \mathrm{D}_{2n}, \mu_1\uparrow\mathrm{D}_{2n}$ with the characters of $\mathrm{D}_{2n}$ in Table \ref{deven} are:
\begin{align*}
&\langle \mu_0 \uparrow \mathrm{D}_{2n}, \chi_1\rangle = 1;
&\langle \mu_1 \uparrow \mathrm{D}_{2n}, \chi_1\rangle = 0;\\
&\langle \mu_0 \uparrow \mathrm{D}_{2n}, \chi_2\rangle = 0;
&\langle \mu_1 \uparrow \mathrm{D}_{2n}, \chi_2\rangle = 1;\\
&\langle \mu_0 \uparrow \mathrm{D}_{2n}, \chi_3\rangle = 1;
&\langle \mu_1 \uparrow \mathrm{D}_{2n}, \chi_3\rangle = 0;\\
&\langle \mu_0 \uparrow \mathrm{D}_{2n}, \chi_4\rangle = 0;
&\langle \mu_1 \uparrow \mathrm{D}_{2n}, \chi_4\rangle = 1;\\
&\langle \mu_0 \uparrow \mathrm{D}_{2n}, \psi_j\rangle = 1;
&\langle \mu_1 \uparrow \mathrm{D}_{2n}, \psi_j\rangle = 1.
\end{align*}
Thus $(\mathrm{D}_{2n}, \langle ba^i\rangle)$ is a strong Gelfand pair for all $i$, independent of the parity of $n$.
\qed

\begin{prop}\label{DinD}
For $n \geq 3$, if $1 \leq m \mid n$, then $\mathrm{D}_{2m} \leq \mathrm{D}_{2n}$, and $(\mathrm{D}_{2n}, \mathrm{D}_{2m})$ is a strong Gelfand pair.
\end{prop}
\noindent\textit{Proof.} By Proposition \ref{Dsmall} for all $i$ we have $\langle ba^i\rangle$ is a strong Gelfand subgroup of $\mathrm{D}_{2n}$ when $n\geq3$. 
Given $\mathrm{D}_{2m}\leq \mathrm{D}_{2n}$, there is an involution $r\in\mathrm{D}_{2m}$ corresponding to an element $ba^i$ in $\mathrm{D}_{2n}$. Thus $(\mathrm{D}_{2n}, \langle r\rangle)$ is a strong Gelfand pair, and by Lemma \ref{stack} $(\mathrm{D}_{2n}, \mathrm{D}_{2m})$ is a strong Gelfand pair.
\qed

\begin{prop}\label{Dcyclic}
Let $n\geq 3$. The pair $(\mathrm{D}_{2n}, \langle a \rangle)$ is a strong Gelfand pair.
\end{prop}
\noindent\textit{Proof. } We know that for all $n$, $\mathcal{C}_n \cong \langle a \rangle \leq \mathrm{D}_{2n}$. Suppose first that $n$ is odd. Then inducing the characters of $\langle a \rangle$ from Table \ref{C} to $\mathrm{D}_{2n}$ gives $\mu_0\uparrow \mathrm{D}_{2n} = \chi_1+\chi_2$ and $\mu_k\uparrow\mathrm{D}_{2n} = \psi_k$. It is then immediate that the inner products are



$$
\langle \mu_0\uparrow \mathrm{D}_{2n}, \chi_1\rangle=1; \quad
\langle \mu_0\uparrow \mathrm{D}_{2n}, \chi_2\rangle=1; \quad
\langle \mu_0\uparrow \mathrm{D}_{2n}, \psi_j\rangle=0;
$$
$$
\langle \mu_k\uparrow \mathrm{D}_{2n}, \chi_1\rangle=0; \quad
\langle \mu_k\uparrow \mathrm{D}_{2n}, \chi_2\rangle=0; \quad
\langle \mu_k\uparrow \mathrm{D}_{2n}, \psi_j\rangle=\begin{cases} 0 \text{ if } j \neq k\\ 1 \text{ if } j=k.\end{cases}
$$

Suppose that $n=2m$. Rather than inducing we will use Frobenius reciprocity to restrict the characters from Table \ref{deven} to $\mathcal{C}_m \cong \langle a^2 \rangle$; see Table \ref{DrestC}. 

\begin{table}[h!]
\centering
\caption{Characters of $\mathrm{D}_{2n}$ restricted to $\mathcal{C}_m\cong\langle a^2\rangle$ for $n$ even}
\label{DrestC}
\vspace{5pt}
\begin{tabular}{c|ccccc}
 & $1$ & $a^2$ & $a^4$  & $\ldots$ & $a^{2(m-1)}$\\
\hline
$\chi_1 \downarrow \langle a^2 \rangle$ & $1$ & $1$ & $1$ & $\ldots$ & $1$\\
$\chi_2 \downarrow \langle a^2 \rangle$ & $1$ & $1$ & $1$ & $\ldots$ & $1$\\
$\chi_3 \downarrow \langle a^2 \rangle$ & $1$ & $1$ & $1$ & $\ldots$ & $1$\\
$\chi_4 \downarrow \langle a^2 \rangle$ & $1$ & $1$ & $1$ & $\ldots$ & $1$\\
$\psi_j \downarrow \langle a^2 \rangle$ & $2$ & $e^{\frac{2\pi i 2j}{n}}+e^{-\frac{2\pi i 2j}{n}}$ & $\ldots$ & $\ldots$ & $e^{\frac{2\pi i 2(m-1)j}{n}}+e^{-\frac{2\pi i 2(m-1)j}{n}}$
\end{tabular}
\end{table}

The inner products of the characters of $\mathrm{D}_{2n}$ restricted to $\mathcal{C}_m$ with the characters of $\mathcal{C}_m$ are:

$$
\langle \chi_1\downarrow \langle a^2\rangle, \mu_0\rangle = 1; \quad \langle \chi_1 \downarrow \langle a^2\rangle, \mu_k\rangle =0; \quad \langle \psi_j\downarrow\langle a^2\rangle, \mu_0\rangle = 0;
$$
$$
\langle \psi_j\downarrow \langle a^2\rangle, \mu_k\rangle = \begin{cases} 1 \text{ if $k=j$ or $k+j=m$}\\ 0 \text{ otherwise.}\end{cases}
$$

We note that since all the linear characters of $\mathrm{D}_{2n}$ restrict to the same character of $\mathcal{C}_m$, we need not calculate the products for all four. Since no inner product is greater than one, the pair $(\mathrm{D}_{2n}, \langle a^2\rangle)$ is a strong Gelfand pair. Then by Lemma \ref{stack}, $(\mathrm{D}_{2n}, \langle a\rangle)$ is a strong Gelfand pair.
\qed

\medskip
Notice we have Corollary \ref{cycsquare} as an immediate consequence of our proof.

\begin{cor}\label{cycsquare}
Let $n\geq 3$. If $n=2m$ then $(\mathrm{D}_{2n}, \langle a^2\rangle)$ is a strong Gelfand pair.
\end{cor}

\begin{prop}\label{dnomore}
Let $n\geq 3$. For $n$ odd, no proper subgroup of $\langle a \rangle \leq \mathrm{D}_{2n}$ is a strong Gelfand subgroup. For $n$ even, no proper subgroup of $\langle a^2\rangle$ is a strong Gelfand subgroup.
\end{prop}
\noindent\textit{Proof. } Suppose that $n$ is odd. Let $m\mid n$ with $1< m\neq n$ and notice that $m < (n-1)/2$ and $n/m=k$. Further $\mathcal{C}_m \leq \mathrm{D}_{2n}$ and so we can restrict the character $\psi_m$ to $\mathcal{C}_m$ to get:

\begin{table}[h!]
\centering
\vspace{5pt}
\begin{tabular}{c|ccccc}
 & $1$ & $a^k$ & $a^{2k}$ & $\ldots$ & $a^{(m-1)k}$\\
\hline
$\psi_m \downarrow \mathcal{C}_{m}$ & $2$ & $2$ & $2$ & $\ldots$ & $2$\\
\end{tabular}
\end{table}
\noindent
So $\langle\psi_m\downarrow\mathcal{C}_m,\mu_0\rangle=2$, which shows that $(\mathrm{D}_{2n}, \mathcal{C}_m)$ is not a strong Gelfand pair.

Suppose $n=2m$ and let $d\mid n$ with $1< d\neq n$. We let $n/d=k$ and notice that $d \leq n/2=m$. Further $\mathcal{C}_d \leq \mathrm{D}_{2n}$ and so if $d\neq m$ we can restrict the character $\psi_d$ to $\mathcal{C}_d$ to obtain: 

\begin{table}[h!]
\centering
\vspace{5pt}
\begin{tabular}{c|ccccc}
 & $1$ & $a^k$ & $a^{2k}$ & $\ldots$ & $a^{(p-1)k}$\\
\hline
$\psi_d \downarrow \mathcal{C}_{d}$ & $2$ & $2$ & $2$ & $\ldots$ & $2$\\
\end{tabular}
\end{table}

Then $\langle\psi_d\downarrow\mathcal{C}_d,\mu_0\rangle=2$, so $(\mathrm{D}_{2n},\mathcal{C}_d)$ is not a strong Gelfand pair when $d\neq m$.
\qed

\medskip
\noindent\textsc{Proof of Theorem \ref{dsgp}.} By Lemma \ref{dsub} we need only consider cyclic and dihedral subgroups. Proposition \ref{Dsmall} covers case (i), Proposition \ref{DinD} covers case (ii), and Proposition \ref{Dcyclic} covers case (iii). Corollary \ref{cycsquare} shows for $n$ even that case (iv) holds. Any remaining subgroup $H$ must be cyclic and a proper subgroup of $\langle a \rangle$. By Proposition \ref{dnomore}, $H$ cannot be a strong Gelfand subgroup.
\qed

\section{Dicyclic Groups}

If $n\geq 3$ is odd, the character table for $\mathrm{Dic}_{4n}$ is given in Table \ref{dicodd}, where $1\leq r \leq n-1$, and $1\leq j , k\leq (n-1)/2$.

\begin{table}[h!]
\centering
\caption{Character Table for $\mathrm{Dic}_{4n}$ with $n$ odd}
\label{dicodd}
\vspace{5pt}
\begin{tabular}{c|ccccc}
& $1$ & $b^2$ & $a^r$ & $b$ & $ba$\\
\hline
$\theta_1$ & $1$ & $1$ & $1$ & $1$ & $1$\\
$\theta_2$ & $1$ & $1$ & $1$ & $-1$ & $-1$\\
$\theta_3$ & $1$ & $-1$ & $(-1)^r$ & $i$ & $-i$\\
$\theta_4$ & $1$ & $-1$ & $(-1)^r$ & $-i$ & $i$\\
$\pi_j$ & $2$ & $2$ & $e^{\frac{2\pi ijr}{n}} + e^{-\frac{2\pi ijr}{n}}$ & $0$ & $0$\\
$\gamma_k$ & $2$ & $-2$ & $e^{\frac{\pi irk}{n}}+e^{-\frac{\pi irk}{n}}$ & $0$ & $0$ 
\end{tabular}
\end{table}

When $n=2m$ the character table of $\mathrm{Dic}_{4n}$ is the same as the character table for $\mathrm{D}_{4n}$; see Table \ref{deven}. It is sufficient to consider $n \geq 2$ since when $n=1$ the group $\mathrm{Dic}_{4}\cong\mathcal{C}_4$ is abelian. 

\begin{prop}\label{dicsmall}
Let $n\geq 2$. The pair $(\mathrm{Dic}_{4n}, \langle ba^i \rangle)$ is a strong Gelfand pair.
\end{prop}
\noindent\textit{Proof.} First suppose that $n\geq 3$ is odd. Then we can restrict the characters from Table \ref{dicodd} to $\langle b \rangle$; see Table \ref{dicdowncodd}. It is sufficient to consider only $b$ since $ba^i$ is conjugate to $b$ if $i$ is even, and if $i$ is odd we see this only interchanges the values of $\theta_1, \theta_3$ with $\theta_2, \theta_4$, respectively.

\begin{table}[h!]
\centering
\caption{Characters of $\mathrm{Dic}_{4n}$ with $n$ odd, restricted to $\mathcal{C}_4\cong\langle b\rangle$}
\label{dicdowncodd}
\vspace{5pt}
\begin{tabular}{c|cccc}
& $1$ & $b^2$ & $b$ & $ba$\\
\hline
$\theta_1\downarrow\mathcal{C}_4$ & $1$ & $1$ & $1$ & $1$\\
$\theta_2\downarrow\mathcal{C}_4$ & $1$ & $1$ & $-1$ & $-1$\\
$\theta_3\downarrow\mathcal{C}_4$ & $1$ & $-1$ & $i$ & $-i$\\
$\theta_4\downarrow\mathcal{C}_4$ & $1$ & $-1$ & $-i$ & $i$\\
$\pi_j\downarrow\mathcal{C}_4$ & $2$ & $2$ & $0$ & $0$\\
$\gamma_k\downarrow\mathcal{C}_4$ & $2$ & $-2$ & $0$ & $0$ 
\end{tabular}
\end{table}

Taking the inner product with the characters of $\langle b \rangle$ from Table \ref{C} yields
\begin{align*}
\langle \theta_1\downarrow\mathcal{C}_4, \mu_0\rangle =1;\quad
\langle \theta_1\downarrow\mathcal{C}_4, \mu_1\rangle =0;\quad
\langle \theta_1\downarrow\mathcal{C}_4, \mu_2\rangle =0;\quad
\langle \theta_1\downarrow\mathcal{C}_4, \mu_3\rangle =0;\\
\langle \theta_2\downarrow\mathcal{C}_4, \mu_0\rangle =0;\quad
\langle \theta_2\downarrow\mathcal{C}_4, \mu_1\rangle =0;\quad
\langle \theta_2\downarrow\mathcal{C}_4, \mu_2\rangle =1;\quad
\langle \theta_2\downarrow\mathcal{C}_4, \mu_3\rangle =0;\\
\langle \theta_3\downarrow\mathcal{C}_4, \mu_0\rangle =0;\quad
\langle \theta_3\downarrow\mathcal{C}_4, \mu_1\rangle =0;\quad
\langle \theta_3\downarrow\mathcal{C}_4, \mu_2\rangle =0;\quad
\langle \theta_3\downarrow\mathcal{C}_4, \mu_3\rangle =1;\\
\langle \theta_4\downarrow\mathcal{C}_4, \mu_0\rangle =0;\quad
\langle \theta_4\downarrow\mathcal{C}_4, \mu_1\rangle =1;\quad
\langle \theta_4\downarrow\mathcal{C}_4, \mu_2\rangle =0;\quad
\langle \theta_4\downarrow\mathcal{C}_4, \mu_3\rangle =0;\\
\langle \pi_j\downarrow\mathcal{C}_4, \mu_0\rangle =1;\quad
\langle \pi_j\downarrow\mathcal{C}_4, \mu_1\rangle =0;\quad
\langle \pi_j\downarrow\mathcal{C}_4, \mu_2\rangle =1;\quad
\langle \pi_j\downarrow\mathcal{C}_4, \mu_3\rangle =0;\\
\langle \gamma_k\downarrow\mathcal{C}_4, \mu_0\rangle =0;\quad
\langle \gamma_k\downarrow\mathcal{C}_4, \mu_1\rangle =1;\quad
\langle \gamma_k\downarrow\mathcal{C}_4, \mu_2\rangle =0;\quad
\langle \gamma_k\downarrow\mathcal{C}_4, \mu_3\rangle =1.
\end{align*}

Thus if $n$ is odd $(\mathrm{Dic}_{4n}, \langle b\rangle)$ is a strong Gelfand pair. Suppose now that $n$ is even. Then restricting the characters of $\mathrm{Dic} _{4n}$ from Table \ref{deven} to $\langle b \rangle$ yields the values in Table \ref{dicdownceven}. It is again sufficient to consider only $b$ in place of $ba^i$ because if $i$ is even they are conjugate, and if $i$ is odd this swaps the values of the characters $\chi_3, \chi_4$. 

\begin{table}[h!]
\centering
\caption{Characters of $\mathrm{Dic}_{4n}$ with $n$ even, restricted to $\mathcal{C}_4\cong\langle b\rangle$}
\label{dicdownceven}
\vspace{5pt}
\begin{tabular}{c|cccc}
& $1$ & $b^2$ & $b$ & $ba$\\
\hline
$\chi_1\downarrow\mathcal{C}_4$ & $1$ & $1$ & $1$ & $0$\\
$\chi_2\downarrow\mathcal{C}_4$ & $1$ & $1$ & $-1$ & $0$\\
$\chi_3\downarrow\mathcal{C}_4$ & $1$ & $1$ & $1$ & $0$\\
$\chi_4\downarrow\mathcal{C}_4$ & $1$ & $1$ & $-1$ & $0$\\
$\psi_j\downarrow\mathcal{C}_4$ & $2$ & $2(-1)^j$ & $0$ & $0$
\end{tabular}
\end{table}

Of note is that $\chi_1\downarrow\mathcal{C}_4=\chi_3\downarrow\mathcal{C}_4$ and $\chi_2\downarrow\mathcal{C}_4=\chi_4\downarrow\mathcal{C}_4$. Armed with these values, we can now calculate that the inner products between characters of $\langle b \rangle$ and characters of $\mathrm{Dic}_{4n}$ restricted to $\langle b\rangle$ are
\begin{align*}
\langle \mu_0, \chi_1\downarrow\langle b \rangle\rangle = 1;\quad
\langle \mu_0, \chi_2\downarrow\langle b \rangle\rangle = 0;\quad
\langle \mu_1, \chi_1\downarrow\langle b \rangle\rangle = 0;\quad
\langle \mu_1, \chi_2\downarrow\langle b \rangle\rangle = 0;\\
\langle \mu_2, \chi_1\downarrow\langle b \rangle\rangle = 0;\quad
\langle \mu_2, \chi_2\downarrow\langle b \rangle\rangle = 1;\quad
\langle \mu_3, \chi_1\downarrow\langle b \rangle\rangle = 0;\quad
\langle \mu_3, \chi_2\downarrow\langle b \rangle\rangle = 0;\\
\langle \mu_0, \psi_j\downarrow\langle b \rangle\rangle = \begin{cases} 1 \text{ if $j$ is even}\\ 0 \text{ if $j$ is odd;} \end{cases}
\langle \mu_1, \psi_j\downarrow\langle b \rangle\rangle = \begin{cases} 0 \text{ if $j$ is even}\\ 1 \text{ if $j$ is odd;} \end{cases}\\
\langle \mu_2, \psi_j\downarrow\langle b \rangle\rangle = \begin{cases} 1 \text{ if $j$ is even}\\ 0 \text{ if $j$ is odd;} \end{cases}
\langle \mu_3, \psi_j\downarrow\langle b \rangle\rangle = \begin{cases} 0 \text{ if $j$ is even}\\ 1 \text{ if $j$ is odd.} \end{cases}
\end{align*}

Thus if $n$ is even then $(\mathrm{Dic}_{4n}, \langle b \rangle)$ is a strong Gelfand pair.

\qed

\begin{prop}\label{dicindic}
Let $n\geq 2$. Any dihedral subgroup $\mathrm{Dic}_{4k}\leq \mathrm{Dic}_{4n}$ is a strong Gelfand subgroup.
\end{prop}
\noindent\textit{Proof. } Follows from Proposition \ref{dicsmall} and Lemma \ref{stack}.\qed

\begin{prop}\label{dicsquare}
Let $n\geq 2$. The pair $(\mathrm{Dic}_{4n}, \langle a^2 \rangle)$ is a strong Gelfand pair.
\end{prop}
\noindent\textit{Proof. } 
Suppose that $n$ is odd. We can restrict the characters of $\mathrm{Dic}_{4n}$ from Table \ref{dicodd} to the subgroup $\langle a^2\rangle$ to get Table \ref{dicoddc} where $1\leq r \leq (n-1)/2$ and $1 \leq  j, k \leq (n-1)/2$.

\begin{table}[h!]
\centering
\caption{Characters of $\mathrm{Dic}_{4n}$ with $n$ odd, restricted to $\mathcal{C}_n\cong\langle a^2\rangle$}
\label{dicoddc}
\vspace{5pt}
\begin{tabular}{c|cc}
& $1$ & $(a^2)^r$\\
\hline
$\theta_1\downarrow\langle a^2 \rangle$ & $1$ & $1$\\
$\theta_2\downarrow\langle a^2 \rangle$ & $1$ & $1$\\
$\theta_3\downarrow\langle a^2 \rangle$ & $1$ & $1$\\
$\theta_4\downarrow\langle a^2 \rangle$ & $1$ & $1$\\
$\pi_j\downarrow\langle a^2 \rangle$ & $2$ & $e^{\frac{4\pi ijr}{n}} + e^{-\frac{4\pi ijr}{n}}$\\
$\gamma_k\downarrow\langle a^2 \rangle$ & $2$ & $e^{\frac{2\pi irk}{n}}+e^{-\frac{2\pi irk}{n}}$
\end{tabular}
\end{table}

We then need to consider fewer products, given that the characters $\theta_1\downarrow\langle a^2\rangle = \theta_2\downarrow\langle a^2\rangle = \theta_i\downarrow\langle a^2\rangle =\theta_3\downarrow\langle a^2\rangle =\theta_4\downarrow\langle a^2\rangle$. The inner products with the characters of Table \ref{C} are
\begin{align*}
\langle\theta_1\downarrow\langle a^2 \rangle, \mu_0\rangle = 1;\quad
\langle\theta_1\downarrow\langle a^2 \rangle, \mu_s\rangle = 0;\quad
\langle\pi_j\downarrow\langle a^2 \rangle, \mu_0\rangle = 0;\quad 
\langle\gamma_k\downarrow\langle a^2 \rangle, \mu_0\rangle = 0;
\end{align*}
\begin{align*}
&\langle\pi_j\downarrow\langle a^2 \rangle, \mu_s\rangle = \begin{cases}1\text{ if } s-2j=0\text{ or } s+2j=n\\0 \text{ otherwise;}\end{cases}\\
&\langle\gamma_k\downarrow\langle a^2 \rangle, \mu_s\rangle = \begin{cases}1 \text{ if } s+k=n \text{ or } s-k =0\\0\text{ otherwise.}\end{cases}
\end{align*}

Thus if $n$ is odd, $(\mathrm{Dic}_{4n}, \langle a^2\rangle)$ is a strong Gelfand pair.
Supposing instead that $n$ is even, we get the restriction of characters from $\mathrm{Dic}_{4n}$ to $\langle a^2\rangle$ shown in Table \ref{dicevenc}. 

\begin{table}[h!]
\centering
\caption{Characters of $\mathrm{Dic}_{4n}$, $n=2m$, restricted to $\mathcal{C}_n\cong\langle a^2\rangle$}
\label{dicevenc}
\vspace{5pt}
\begin{tabular}{c|ccc}
& $1$ & $a^{2n}$ & $(a^2)^r$\\
\hline
$\chi_1\downarrow\langle a^2 \rangle$ & $1$ & $1$ & $1$\\
$\chi_2\downarrow\langle a^2 \rangle$ & $1$ & $1$ & $1$\\
$\chi_3\downarrow\langle a^2 \rangle$ & $1$ & $1$ & $1$\\
$\chi_4\downarrow\langle a^2 \rangle$ & $1$ & $1$ & $1$\\
$\psi_j\downarrow \langle a^2\rangle$ & $2$ & $2(-1)^j$ & $e^{\frac{\pi i r j}{n}} + e^{-\frac{\pi r i j}{n}}$
\end{tabular}
\end{table}

Since $\chi_1\downarrow \langle a^2\rangle = \chi_2\downarrow \langle a^2\rangle = \chi_3\downarrow \langle a^2\rangle = \chi_4\downarrow \langle a^2\rangle$ it is sufficient to check the inner products of these restricted characters with the characters of $\langle a^2\rangle$ given in Table \ref{C}:
\begin{align*}
\langle \chi_1\downarrow \langle a^2\rangle, \mu_0\rangle = 1;\quad
\langle \chi_1\downarrow \langle a^2\rangle, \mu_k\rangle = 0;\quad
\langle \psi_j\downarrow \langle a^2\rangle, \mu_0\rangle = 0;\\
\langle \psi_j\downarrow \langle a^2\rangle, \mu_k\rangle = \begin{cases}1 \text{ if } k=j \text{ or } k+j=m\\ 0 \text{ otherwise.}\end{cases}
\end{align*}

Thus $(\mathrm{Dic}_{4n}, \langle a^2\rangle)$ is a strong Gelfand pair when $n$ is even. \qed

\begin{cor}\label{diccyclic}
Let $n\geq 2$. The pair $(\mathrm{Dic}_{4n}, \langle a\rangle)$ is a strong Gelfand pair.
\end{cor}
\noindent\textit{Proof.} Immediate from Proposition \ref{dicsquare} and Lemma \ref{stack}. \qed

\begin{prop}\label{dicnomore}
Let $n\geq 2$. No proper subgroup of $\langle a^2\rangle \leq \mathrm{Dic}_{4n}$ is a strong Gelfand subgroup.
\end{prop}
\noindent\textit{Proof.}
Suppose $n$ is odd and take $d\in\mathbb{Z}^+$ such that $1 < d$ and $d\mid 2n$. Then $\langle a^{2n/d}\rangle \cong \mathcal{C}_d \lneq \mathcal{C}_{2n} \cong \langle a\rangle \leq \mathrm{Dic}_{4n}$. Observe that $\langle \pi_d\downarrow \mathcal{C}_d\rangle = 2$. So $(\mathrm{Dic}_{4n}, \langle a^{2n/d}\rangle)$ is not a strong Gelfand pair. 

Suppose now that $n$ is even, and once again take $d\in\mathbb{Z}^+$ such that $1 < d$ and $d\mid n$. Then $\langle a^{2n/d}\rangle \cong \mathcal{C}_d \lneq \mathcal{C}_{2n} \cong \langle a\rangle \leq \mathrm{Dic}_{4n}$. So restricting the character $\psi_d$ to $\langle a^{2n/d}\rangle$ gives us $\langle \psi_d\downarrow \langle a^{2n/d}\rangle, \mu_o\rangle=2$, showing that $(\mathrm{Dic}_{4n}, \langle a^{2n/d}\rangle)$ is not a strong Gelfand pair.
\qed

\begin{cor}\label{diccenter}
Let $n\geq 2$. The pair $(\mathrm{Dic}_{4n}, \langle b^2\rangle)$ is not a strong Gelfand pair. 
\end{cor}
\noindent\textit{Proof.} Since $a^n=b^2$, this follows immediately from Corollary \ref{diccyclic}.\qed

\medskip
\noindent\textsc{Proof of Theorem \ref{dicsgp}.} Proposition \ref{dicsmall} covers case (i), Proposition \ref{dicindic} covers case (ii), and Proposition \ref{dicsquare} with Corollary \ref{diccyclic} covers case (iii). Let $H\leq \mathrm{Dic}_{4n}$ be a subgroup not in cases (i), (ii), and (iii). By Lemma \ref{dicsub} any remaining subgroup of a dicyclic group is cyclic or dihedral, so $H$ is cyclic. Then $H$ is generated by a power of $a$ or a power of $ba^i$. If $H$ is generated by a power of $a$, then Proposition \ref{dicnomore} says $(\mathrm{Dic}_{4n}, H)$ is not a strong Gelfand pair. If $H$ is generated by a power of $ba^i$ then $H$ is either as in case (i), or $H = \langle b^2\rangle$, giving by Corollary \ref{diccenter} that $(\mathrm{Dic}_{4n}, H)$ is not a strong Gelfand pair. \qed

\end{document}